\documentclass[11pt]{article}
\usepackage{amssymb}
\usepackage[pdftex]{graphicx}
\usepackage{epsfig}
\usepackage{cite}
\newcommand{\AmS}{{\protect\the\textfont2
  A\kern-.1667em\lower.5ex\hbox{M}\kern-.125emS}}
\newcommand{\ba}{\begin{array}}
\newcommand{\ea}{\end{array}}

\def\beq{\begin{equation}}   
\def\eeq{\end{equation}}

\def\bea{\begin{eqnarray}}
\def\eea{\end{eqnarray}}

\def\beq{\begin{equation}}   
\def\eeq{\end{equation}}

\def\bea{\begin{eqnarray}}
\def\eea{\end{eqnarray}}

\begin{document}
\begin{titlepage}

\begin{flushright}
CERN-TH-2018-083\\
\end{flushright}

\begin{center}
\vspace{2.7cm}
{\Large{\bf A heuristic study of the distribution of primes}}
\end{center}
\begin{center}
{\Large{\bf in short and not-so-short intervals}}
\end{center}

\vspace{1cm}

\begin{center}

{\bf Miguel-Angel 
Sanchis-Lozano$^{\rm a,b,}$\footnote{Email 
address: Miguel.Angel.Sanchis@ific.uv.es}}
\vspace{1.5cm}\\
\it 
{$^{\rm a}$ Theoretical Physics Department, CERN, 1211 Geneva 23, Switzerland\\
$^{\rm b}$ IFIC, Centro Mixto CSIC-Universitat de Val\`encia, Dr. Moliner 50, 46100 Burjassot, Spain}

\end{center}

\vspace{0.5cm}

\begin{abstract}

A numerical study on the distributions of primes
in short intervals of length $h$
over the natural numbers $N$ is presented. Based on Cram\'er's model in Number Theory, 
we obtain a heuristic expression applicable when $h \gg \log{N}$ but $h \ll N$, providing
support to the Montgomery and Soundararajan conjecture on the variance
of the prime distribution at this scale.
\end{abstract}

\begin{center}

\today

\end{center}

\end{titlepage}


\begin{flushright}
{\it Mathematics and physics, two communities divided by a common language}\\
Modified version of a quotation attributed to George Bernard Shaw
\end{flushright}
\vskip 0.1cm

\section{Introduction}
The study of the appearance of
primes over the natural numbers is a topic of the greatest importance in
Number Theory \cite{Har08}. Historically Gauss ans Legendre conjectured independently, based on
pure empirical evidence,  that the number of
primes below a given integer $x$, denoted as $\pi(x)$, follows a behaviour 
somewhat like $x/\log{x}$.
Later, Hadamard and
de la Vall\'ee Poussin independently proved that
\beq
\pi(x)=\frac{x}{\log{x}}\ \ (x \to \infty)
\eeq
known as the Prime Number Theorem (see \cite{Gran95} and references therein). 
An alternative and more precise version of it states
that, for large $x$,  $\pi(x)$ is asymptotic to the logarithm integral function ${\rm Li}(x)=\int_2^xdx/\log{x}$. 
As is well known, the correcting
terms to this estimate are amazingly related to the famous hypothesis on the non-trivial zeros of the Riemann function,
constituting a long-standing hot topic in mathematics.

Needless to say, not only expected values of prime numbers in intervals are of interest, but also 
their fluctuations about mean values.  In particular, the understanding of the 
distribution of primes in short (and somewhat longer) intervals remains an important
issue in the theory of prime numbers.

\subsection{Cram\'er's model}

It might appear somewhat surprising to talk about ``models'' in mathematics, as 
models are commonly associated to the (only approximate alas!) description of the physical world. Nonetheless, 
the complexity of certain mathematical problems, particularly in Number Theory, would require the
use of approximations which are expected to become, e.g., asymtotically accurate. In fact, models 
for random primes can be
very effective to confidently give an answer to a long list of open questions, whilst not possessing
a clear way forward to rigorously confirm these answers.  

Such models are based on taking some statistical distribution
and replacing it by a {\em model} distribution that is easier to compute with. 
For example, in the Cram\'er model, prime numbers behave like independent random variables with
probability $q=1/\log{x}$ if $x$ is prime, and with probability
$1-1/\log{x}$ if $x$ is composite, which somewhat reminds towing a biased coin
in a series of measurements. 
Of course, the reader may object that the 
probability that $x$ and $x+1$ are both primes must be zero, while the Cram\'er model
assigns this event a finite probability. Taking into account this caveat, one 
should think of rather large numbers 
in order to match Cram\'er's predictions and numerical results, as it happens remarkably well.

Furthermore, letting the number of trials $n$ tend to $\infty$, while keeping
$nq$ fixed, one gets a Poisson distribution. Thus, for any fixed real $\lambda >0$, and integer $k \ge 0$, 
Cram\'er's model implies
\beq
\#\{{\rm integers}\ x \leq N:\ \pi(x+\lambda\log{x})-\pi(x)=k\} \sim N\lambda^ke^{-\lambda}/k!\ \ \ (N \to \infty)
\eeq
In other words, the interval $(x,x+\lambda\log{x})$ contains, on average, $\lambda$ primes as
it represents the mean of a Poissonian distribution of primes.   

Gallagher \cite{Gal76} conditionally proved this result by studying the moments
of the prime distribution in short intervals of length $h$, 
assuming the prime $k$-tuple conjecture by Hardy and Littlewood \cite{HL23}. On the other hand, 
nothing prevents in the Cram\'er model that $h>\log{N}$ with $h \ll N$, i.e. somewhat larger intervals. 
At this scale, Goldston and
Montgomery first \cite{GM84}, and Montgomery and Soundararajan later \cite{MS04},
conjectured that Cram\'er's model overestimates the variance by a factor 
$\sim\ \log{N}/\log{(N/h)}$.

\newpage

\begin{figure}[t!]
\begin{center}
\includegraphics[scale=0.3]{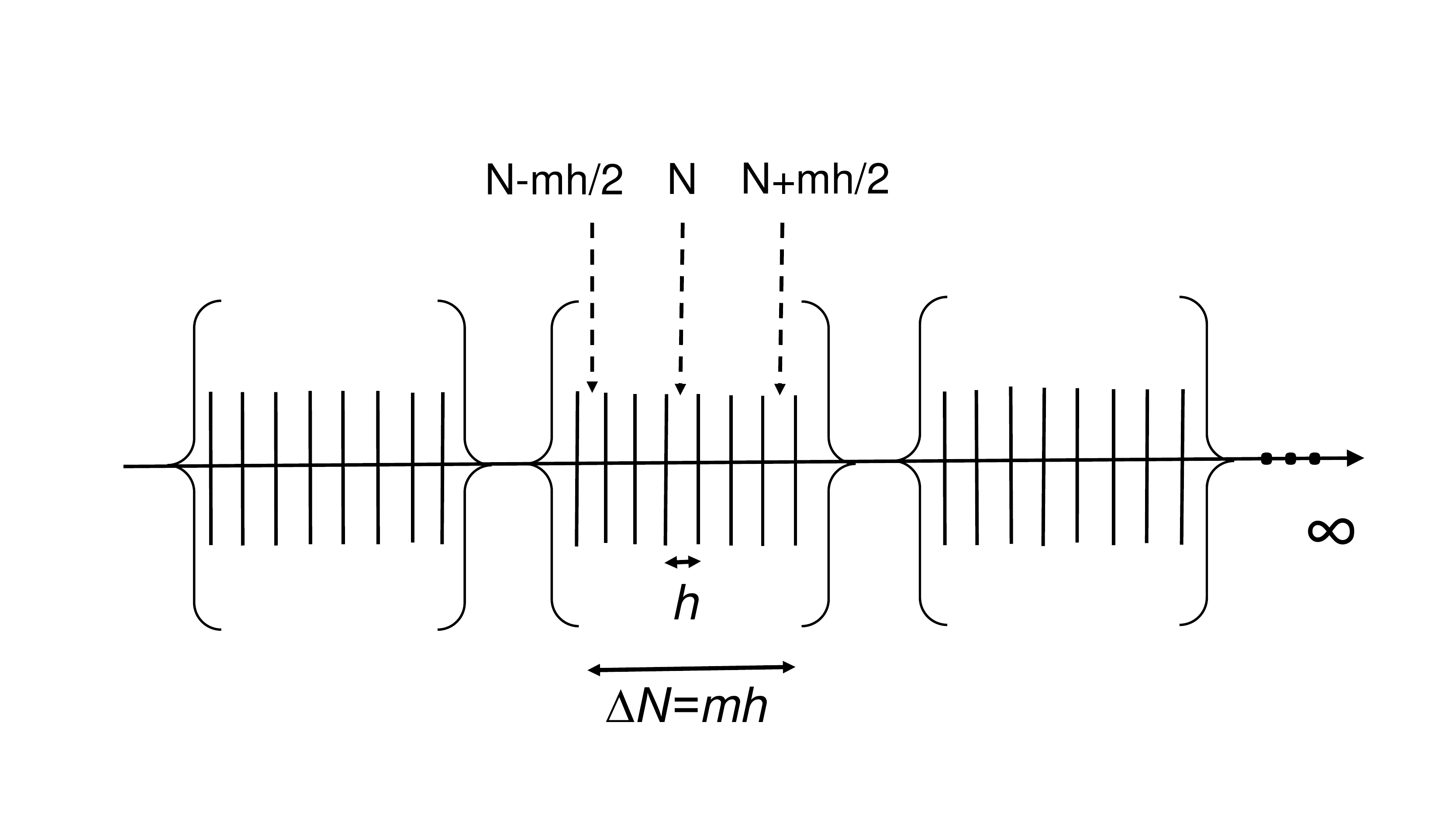}
\caption{Schematic view of sets (in brackets) of $m$ intervals of {\em fixed} length $h$. 
In every set we numerically compute the mean and
variance to be associated to $N$ at the center of the set, for simplicity. In our notation,  
a sample is formed by ensembles of such sets of different lengths $h$, for a given $m$.}
\label{fig:setup}
\end{center}
\end{figure}

\section{Empirical procedure}
In this work, we deal with sets of $m$ intervals of fixed length $h$ centered about a certain value of $N$. 
Although prime numbers are deterministic, for some purposes they
can be viewed as pseudo-random numbers. Hence we look upon
the numbers of primes in intervals as empirical data
coming out from a kind of ``counting experiment''.  Thus we numerically compute the mean $\langle p \rangle$
and the variance $\sigma_p^ 2$ of the primes corresponding to each set of intervals (see Fig.1). In order
to assess the accuracy of our computations, we will estimate
statistical fluctuations and systematic uncertainties (``errors'', as usually known in physics),
as later discussed.

By letting $N$ vary over the natural numbers keeping $h$ and $m$ fixed, we
obtain a series of values for the mean $\langle p \rangle$ and the 
variance $\sigma_p^ 2$, thus computing the normalized variance as the ratio:
\beq
w=\frac{\sigma_p^ 2}{\langle p \rangle}=\ \frac{\langle p^2 \rangle-\langle p \rangle^2}{\langle p \rangle}
\eeq 
for different values of $N$.

Expecting a Poissonian behaviour for $N \to \infty$, 
we parametrize $w$ as
\beq
w\ =\ 1-\frac{b(h,m)}{\log{N}}\ ,
\eeq
where $b(h,m)$ depends on both the interval length $h$ and the number $m$ of intervals 
of each set.
 
Since in our statistical analysis we keep the interval length fixed for each set
about $N$, the corresponding parameter of the Poisson-like distribution 
changes for different choices of $N$ according to
\beq
\lambda(N) = \frac{h}{\log{N}}\ ,
\eeq
representing the expected mean associated to the whole set of $m$ intervals about $N$, within
systematic uncertainties.
Notice that if $h=1$ obviously ${\langle p^2 \rangle}={\langle p \rangle}$ leading to
\beq\label{eq:w1}
w\ =\ 1\ -\ \frac{1}{\log{N}}\ ,
\eeq
as predicted by the Cram\'er model leading to a binomial distribution with $q=1/\log{N}$, and explicitly
checked here by numerical computation.

\begin{figure}[t!]
\begin{center}
\includegraphics[scale=0.81]{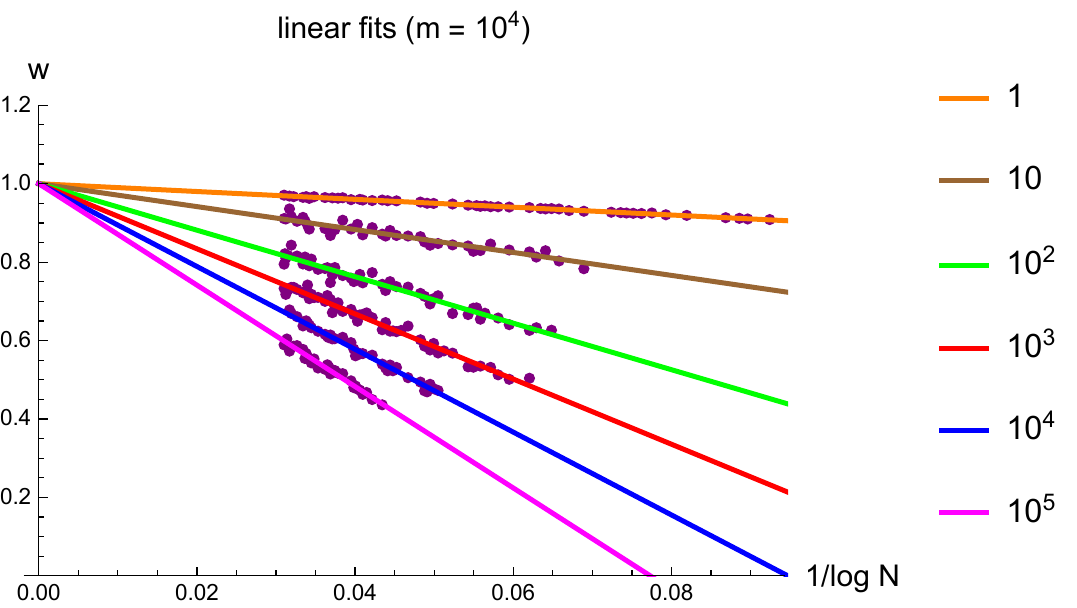}
\hspace{0.5cm}
\includegraphics[scale=0.58]{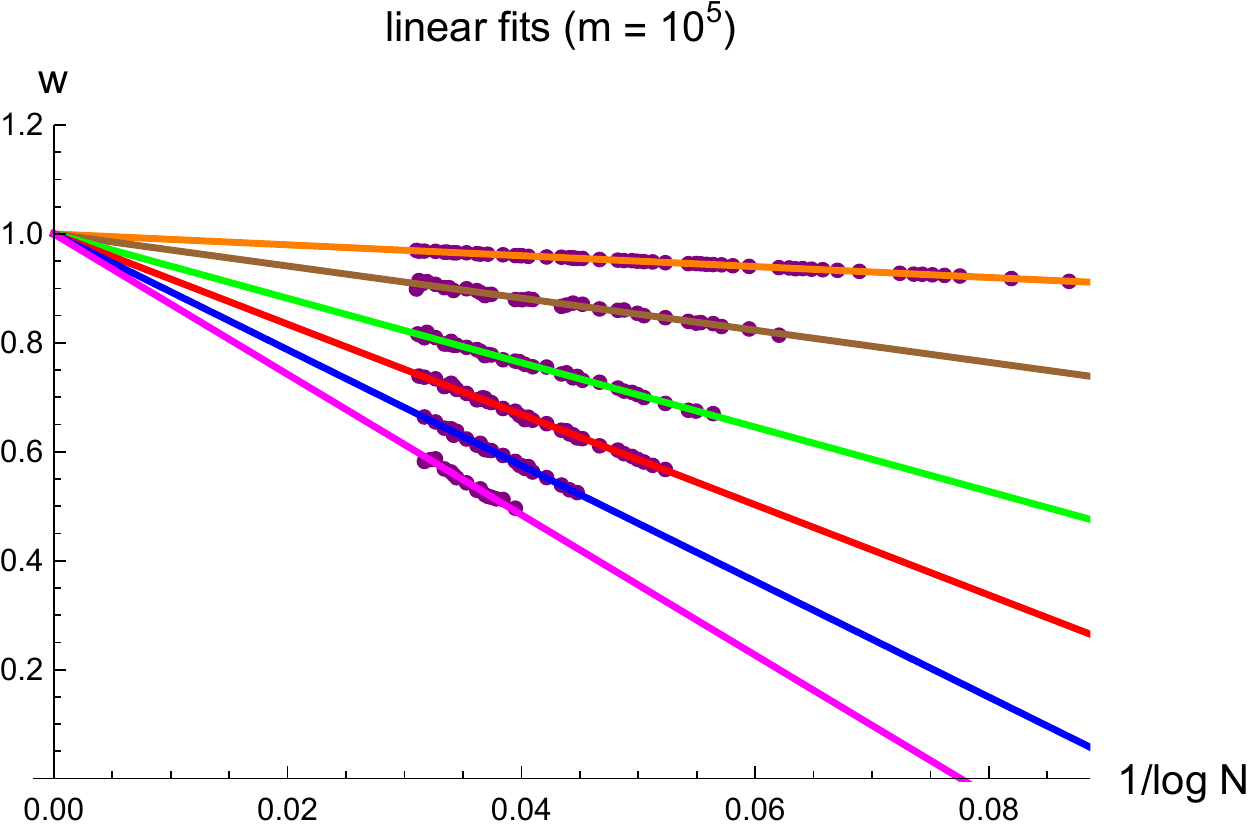}
\caption{Fits of the normalized variance $w=1-b(h,m)/\log{N}$ 
with $h=1,10,10^2,10^3,10^4,10^5$ 
(downwards) for $m=10^4$ (left) and $m=10^5$ (right). Each fit corresponds to an ensemble of sets as schematically shown in Fig.1.}
\label{fig:fit104}
\end{center}
\end{figure}

\subsection{Samples with different numbers of intervals}

In this work, we study three samples with different numbers of intervals: 

\begin{itemize}

\item{sample I:} $m=2 \times 10^3$,
\item{sample II:} $m=10^4$,
\item{sample III:} $m=10^5$.
\end{itemize}

In our fits, the following interval lengths were employed: 
\[
h=1,\ 200,\ 500,\ 1000,\ 2000,\ 5000,\ 10^4,\ 2 \times 10^4,\ 5 \times 10^4,\ 10^5,\ 
2 \times 10^5,\ 5 \times 10^5
\] 

Other interval lengths, namely $h \in \{1,100\}$, were also considered 
in our analysis although they
fall out of the scale of our interest, as later discussed.

\subsection{Error estimates}

As we are thinking or primes as the outcome of a series of ``measurements'', we 
examine here
the sources of systematic and statistical uncertainties in our fits.  
First we find a kind of ``systematic error'' 
when assigning $\langle p \rangle$ to a common value of $N$, while the set of $m$ intervals 
actually spreads over a finite length of $\Delta N=m\times h$ about $N$. To this end, we 
write
\beq
\Delta \lambda = \frac{h}{\log{(N-\Delta N/2)}}-\frac{h}{\log{(N+\Delta N/2)}}=
\frac{2h\ \tanh{}^{-1}(\Delta N/2N)}{\log{(N+\Delta N/2)}\log{(N-\Delta N/2)}}
\eeq
Keeping only first order terms in a Taylor expansion of the above expression for small $\Delta N$
as compared to $N$, we get
\beq
\Delta \lambda \simeq \frac{mh^2}{N(\log{N})^2}\ .
\eeq

The corresponding relative {\em systematic error} can be estimated as (see
Eq.(5))
\beq
\epsilon^{(r)}_{{\rm sys}}\ \simeq \ \frac{\Delta \lambda}{\lambda}\ \simeq\ \frac{mh}{N\log{N}}\ .
\eeq

\begin{table*}[hbt]
\setlength{\tabcolsep}{0.6pc}
\caption{Values of the $\alpha_i$ ($i=1,2$) parameters from the hyperbolic fit 
for the three samples.} \label{tab:Tab1}

\begin{center}
\begin{tabular}{cccc}
  \hline Sample: &  I ($m=2 \times 10^3$)  & II ($m=10^4$) & III  ($m=10^5$) \\
\hline  $\alpha_1$ & $1.00077\pm 0.01111$  & $1.00268 \pm 0.00345$ & $1.00001 \pm 0.00228$ \\
\hline
\hline  $\alpha_2$ & $0.554\pm 0.087$  & $0.614 \pm 0.027$ & $0.5827\pm 0.0179$ \\
\hline

\end{tabular}
\end{center}
\end{table*}

In order to keep this systematic error small enough, the
condition $mh \ll N \log{N}$ has to be satisfied in the whole sample.

On the other hand, considered as ``data'', 
fluctuations of the number of primes around the mean are expected 
for short intervals. From a Poissonian behaviour, 
one expects a relative {\em statistical error} as
\beq
\epsilon^{(r)}_{{\rm stat}}\ \simeq \ \sqrt{\frac{\log{N}}{mh}}
\eeq

Note that errors corresponding to the mean squared and square root mean, used to compute $w$,  
should be (at least) twice the errors of the mean shown above.
In our study, systematic errors of individual points in our fits 
are typically $\lesssim 0.1\%$ (decreasing at large $N$) while statistical errors 
are of order $\lesssim 1\%$ (decreasing for large $m$). 

Extrapolation to higher values of $N$ requires, in turn, larger values of $h$, in order to keep 
both errors under control, i.e. $(\log{N})/h \to 0$  when $N \to \infty$ with $\Delta N \ll N$. These
requirements are consistent with the scale of intervals under study, 
as commented below.

\subsection{Scales}

First of all, we should wonder about the meaning of 
``short'' intervals, and ``not so short'' intervals. In real physics 
experiments, nature usually provides us with 
a scale to compare with the corresponding physical quantity, determining
what one can call ``large'' and what ``small''.  

In our study about primes in intervals, one can intuitively compare $h$ and $N$, 
so that if $h \ll N$, one can talk about short intervals. Nevertheless,
another scales (e.g. taking logarithms) still may be relevant in the problem under consideration. 
For example $\log{N}$ can interpreted as a ``length'' (the inverse of the prime density), somewhat providing a
geometrical sense to the comparison between $h$ and $\log{N}$. 

Thereby, according to \cite{SO04}, let us distinguish the following {\em scales}:
\begin{itemize}
\item {\em Microscopic} scale: when $h$ and $\log{N}$ are not so different, i.e. $h \asymp \log{N}$ 
\footnote{The symbol $ f \asymp g$ means that two real positive  numbers, $C$ and $D$, can be found such that $C <|f|/|g| <D$, in the domain of the functions $f$ and $g$.} which
means that $h$ and $\log{N}$ become asymptotically comparable though $h \ll N$.
\item {\em Mesoscopic} scale: when $\log{h}$ and $\log{N}$ are not so different, i.e. $\log{h} \asymp \log{N}$. This is the case when $h/\log{N} \to \infty$ while $h/N \to 0$ when $N \to \infty$. 
\item {\em Macroscopic} scale: when $h \gg N$. Then no specific distribution law is expected.
\end{itemize}
  
The values of $h$ and $N$ can be considered as inter-related following to the
above-mentioned scales. In our study, mainly focusing on the mesoscopic scale 
where the Montgomery and Soundararajan (MS) conjecture applies, 
we thus require that $h(N)$ always satisfies the condition  $\log{h} \asymp \log{N}$. 
This can be achieved, in our particular case with $h \in \{10^2,10^5\}$ and $N \in \{10^7,10^{14}\}$, by choosing, e.g.,  
$C=\log{10^2}/\log{10^{14}} \simeq 0.1$ and $D=\log{10^5}/\log{10^{7}} \simeq 1$.  
Analogously, $h$ always remains much smaller than $N$, 
as can be easily seen.

\begin{figure}[t!]
\begin{center}
\includegraphics[scale=0.55]{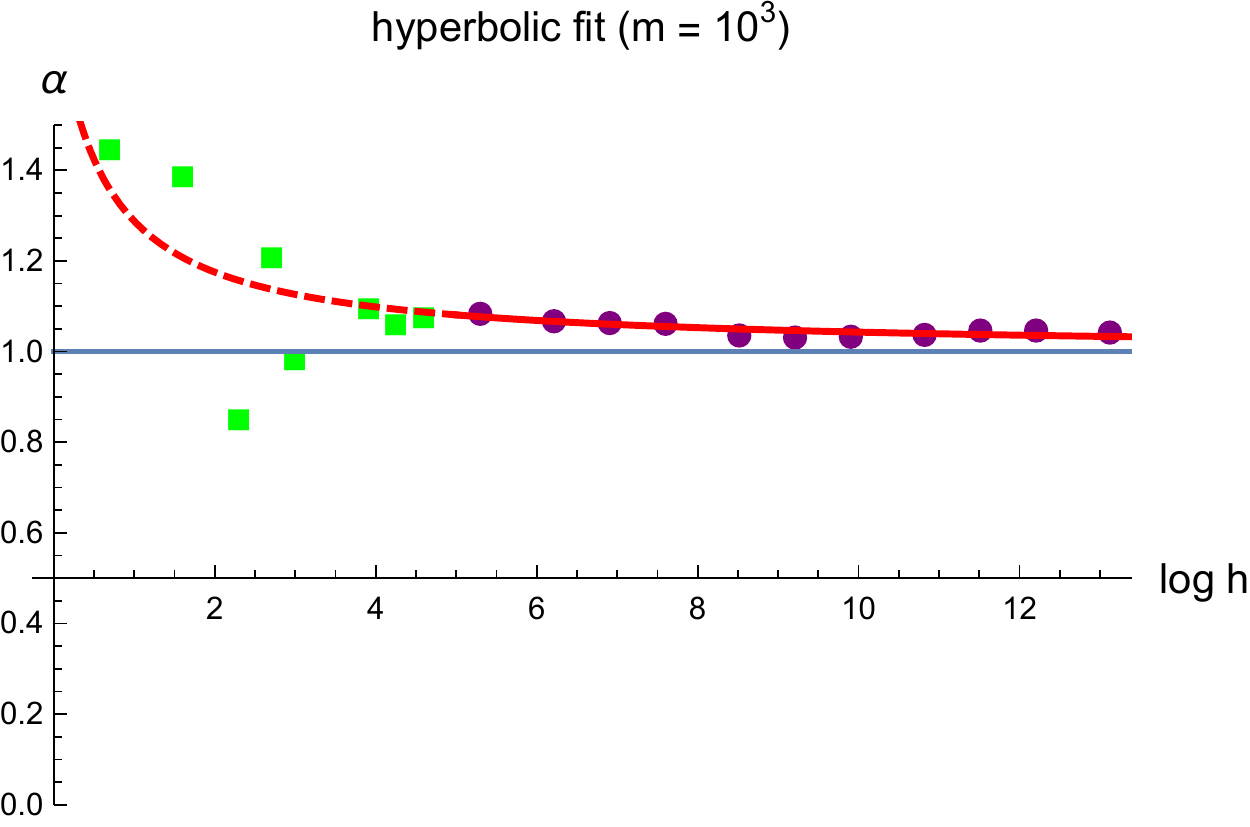}
\hspace{0.5cm}
\includegraphics[scale=0.55]{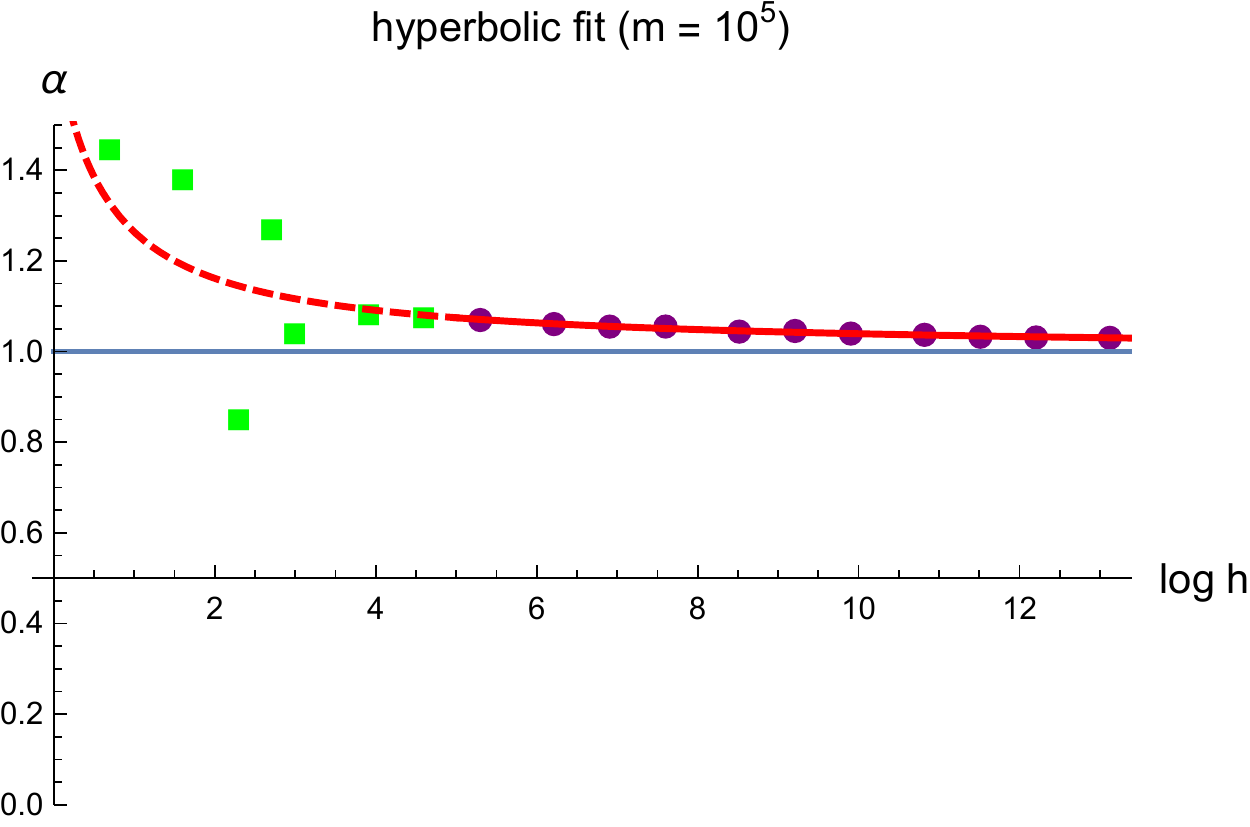}
\caption{Asymptotic limit from sample I (left) and III (right)
of the parameter $\alpha(h)$ defined in Eq.(11) 
as a function of $\log{h}$, using the hyperbolic fit of Eq.(\ref{eq:alphah}). 
The size of the points
are of the order of their statistical error bars. 
Only intervals such that $h > 100$ 
were used in the fits, represented
by dark (magenta) circles. 
Intervals with $h \leq 100$ are also represented as light (green) squares 
but were not used in the fits as they 
do {\em not} satisfy the condition $h \gg \log{N}$.}
\label{fig:contour}
\end{center}
\end{figure}

\section{Results}

In Fig.2 we show the linear plots of $w$ as a function of
$1/\log{N}$ ($w=1-b/\log{N}$), corresponding to the intervals: $h=1,10,10^2,10^3,10^4,10^5$, for
samples I, II and III. Previously, we checked from linear fits $w=a-b/\log{N}$ that
$a$ can be safely set equal to unity within errors in all cases, as expected from the Poisson limit at large $N$.
On the other hand, one can easily see from the figure that statistical flcutuations become smaller 
as $m$ increases.

Having fixed the number $m$ of intervals for each sample, we parametrize 
the dependence of $b(h,m)$ on the interval size $h$ as,
\beq\label{eq:bh}
b(h)=1 +\ \alpha(h)\log{h}
\eeq
as suggested by Eq.(\ref{eq:w1}) and in view of the plots of Fig.2. This fact plays a fundamental role
for later evidence in favour of the MS conjecture.

Next, according to Eqs.(4) and (11) we write
\beq\label{eq:general}
w=1-\frac{1+\alpha\log{h}}{\log{N}}=\frac{\log{(N/h^{\alpha})}}{\log{N}}-\frac{1}{\log{N}}
\eeq
We can study the behaviour of $\alpha(h)$ using the hyperbolic expression 
\beq\label{eq:alphah}
\alpha(h)=\frac{1+\alpha_1\log{h}}{\alpha_2+\log{h}} 
\eeq
for different values of $h$.

The results from the fit for the parameters $\alpha_1$ and $\alpha_2$ 
can be found in Table I
for samples I, II and III. 
 Notice that $\alpha \to 1$ asymptotically in the limit $N \to \infty$,
for all three samples. In fact, $\alpha$
remains quite close to unity for $h \gtrsim 10^2$. 

The accuracy of the fit improves for larger $m$ since it is
mainly determined by statistical errors as discussed before. In Fig.3
we show the curves from the fits of samples I and III. Two regions can be
distinguished : one region corresponds to intervals of length
$h \gtrsim 10^2$ shown by dark circles, used in the fit providing the results of Table I. 
Another region corresponds to intervals
of length $h \lesssim 10^2$ shown by light squares, displaying
broad oscillations and not following the smooth behaviour extrapolated 
from larger $h$ values. Indeed, the latter points
do not satisfy the condition $h \gg \log{N}$, thereby falling out of 
the scale of interest of this paper. 

Now, Taylor expanding Eq.(\ref{eq:alphah}) for small $1/\log{h}$, one gets:
\beq\label{eq:parab}
\alpha(h) \simeq \frac{1-\alpha_1\alpha_2+\alpha_1\log{h}}{\log{h}}\ +\ {\cal O}\biggl(\frac{1}{\log{}^2h}\biggr)
\eeq
Substituting in Eq.(\ref{eq:bh}), one obtains
\beq\label{eq:bh2}
b(h)=2-\alpha_1\alpha_2+\alpha_1\log{h}+\ {\rm smaller\ terms}
\eeq

In case we assume $\alpha_1=1$ and identify $2-\alpha_1\alpha_2=1.417 \pm 0.018$ with
$-B=\gamma_E+\log{(2\pi)}-1=1.414509...$, where $\gamma_E$ is 
the Euler's constant, one can write
\beq
w\ \simeq\ 1-\frac{\log{h}-B}{\log{N}}
\eeq

Therefore, our numerical results (at the given accuracy) support
the MS conjecture stating that the variance of the prime distribution
exhibits the following behaviour: \footnote{The MS conjecture is formulated
in Refs. \cite{MS04} and \cite{SO04} by studying 
the moments of $\psi(x+h)-\psi(x)-h$, where
$\psi(x)=\sum_{n \leq x}\Lambda(n)$ (with
$\Lambda(n)$ denoting the von Mangoldt function) 
instead of $\pi(x+h)-\pi(x)$, therefore not exactly matching
our expression.}
\beq\label{eq:sigmap}
\sigma_p^2\ \sim\ \frac{h}{(\log{N})^2}\ \biggl[\log{\biggl(\frac{N}{h}\biggr)}+B\biggr]\ .
\eeq
for large $h \gg \log{N}$ but $h \ll N$.

Eq.(\ref{eq:sigmap}) reproduces the 
logarithmic dependence, $\log{(N/h)}$,
as predicted by the MS conjecture, for the variance of the 
distribution of prime numbers at the mesoscopic scale. Notice that 
throughout our numerical analysis,  
the condition  $h >\log{N}$ but $h \ll N$.
is fulfilled as $h$ varies typically from $10^2$ up to $10^5$
whereas $\log{N}$ remains of order ${\cal O}(10)$ 
along the $N$
region scanned over the natural numbers.

Let us finally remark that the dependence of the variance on $\log{(N/h)}$   
ultimately comes from the dependence on $\log{h}$ of the slopes $b(h,m)$ of the
fits shown in Fig.2.

\begin{figure}[t!]
\begin{center}
\includegraphics[scale=0.69]{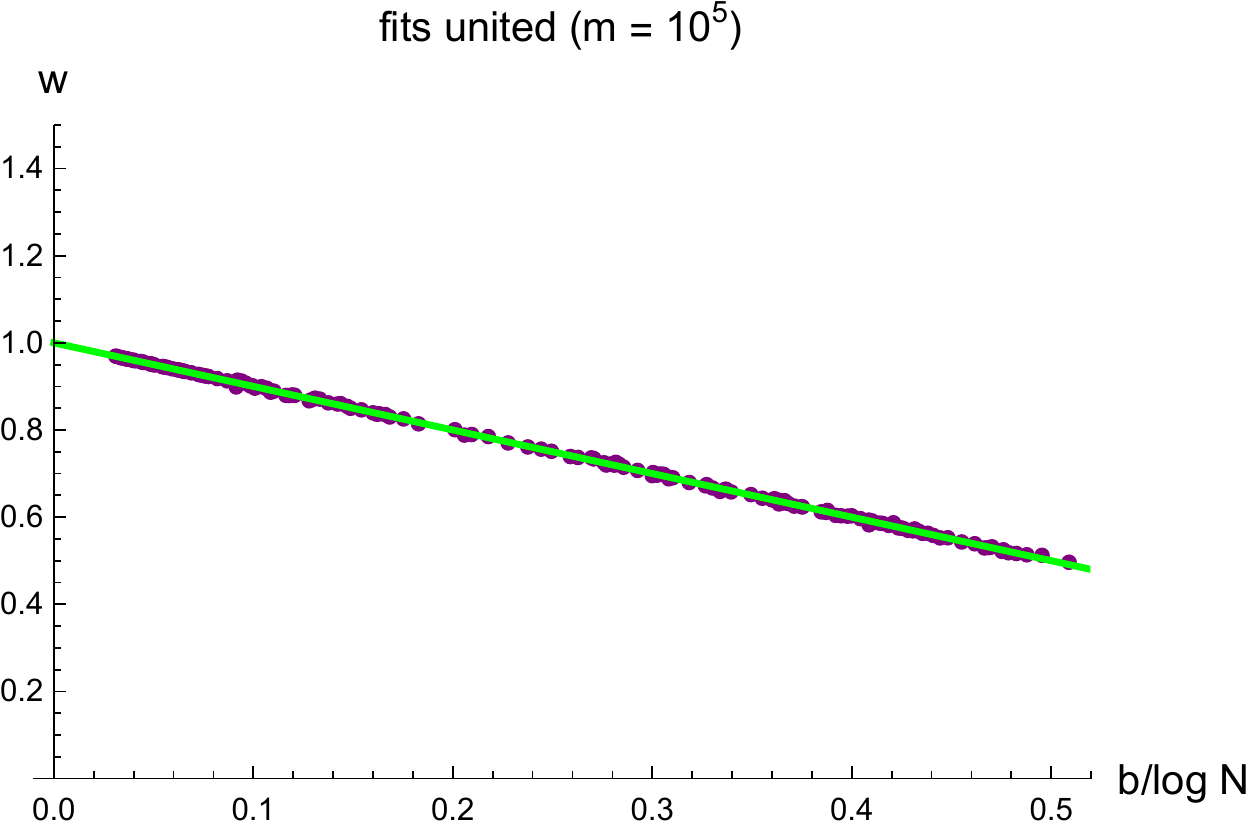}
\caption{Linear fit of the normalized variance $w$ as a function
of $b(h)/\log{N}$ for all sets together from sample III shown in Fig.2 (right).
The result of the fit $w=A-B(b/\log{N})$ yields 
$A=1.000065\pm 0.000410$ and $B=1.00015 \pm 0.00145$.
We look upon this plot as an overall confirmation of 
the consistency of our fits pointing at
a Poisson distribution ($w \to 1$) as $N \to \infty$.}
\label{fig:fit104}
\end{center}
\end{figure}

\newpage

\section{Conclusions}

In this paper, we present a numerical study of the distribution of primes
in short intervals of length $h$ such that $h \gtrsim \log{N}$ 
and $h \ll N$, up to
$N =10^{14}$. In the literature, attempts to check numerically the
Montgomery and Soundararajan conjecture can be found, e.g. in \cite{MS03}, but this specific 
kind of study has not been yet done, to our knowledge.

In our numerical approach we rely on the {\em Mathematica} package
to compute the number of primes in intervals. Let us point out that 
a number of intervals (with different lengths) of  order $10^8$ were generated,  
from which a total of about one thousand points (corresponding to the normalized variance 
$w$ versus $1/\log{N}$) were obtained and used in our fits (see Fig.2, where only a part of them are 
shown). To check the overall consistency of our analysis, in Fig.4 we plot 
$w$ as a function of $b(h)/\log{N}$ from all sets of sample III, 
impressively showing 
that indeed $w \to 1$ as $N \to \infty$, together with the consistency of our results. 

In sum, using the parametrization of Eqs.(\ref{eq:alphah})-(\ref{eq:bh2}),  
we obtain the heuristic expression for the variance of the prime distribution: 
\beq\label{eq:last}
\sigma_p^2 \sim \frac{h}{(\log{N})^2}\ \biggl[\log{\biggl(\frac{N}{h}\biggr)}+B\biggr]
\eeq
for large $h >\log{N}$ but $h \ll N$, representing 
our main result empirically supporting the MS conjecture 
basically implying a dependence 
of the variance on $\log{(N/h)}$. 
Note that Eq.(\ref{eq:last}) does not apply to intervals whose length 
$h \lesssim 100$, since then the parameter $\alpha(h,m)$ in Eq.(\ref{eq:alphah})
fluctuates broadly, as can be seen in Fig.3 (green squares), belonging
to a different scale. Nevertheless, let us remark that
Eq.(\ref{eq:general}) 
provides an interpolating expression between different scales.

In conclusion, even with the relatively small upper values 
of natural numbers reached in our study ($h=10^5$, $N=10^{14}$), one can tentatively
conclude that there is a clear empirical evicence in favour of the 
MS conjecture at the mesoscopic scale.  
Higher values of $h$ in order to check further this conjecture need
higher values of $N$ satisfying the condition $\log{h} \asymp \log{N}$, 
requiring a larger computer capacity.

\subsubsection*{Acknowledgements}

This work has been partially supported by the Spanish MINECO under grants FPA2014-54459-P and
FPA2017-84543-P, by the Severo Ochoa Excellence Program under grant SEV-2014-0398 and by
the Generalitat Valenciana under grant GVPROMETEOII 2014-049.
M.A.S.L. thanks the CERN Theoretical Physics Department, where this work 
was finished, for its warm hospitality.

\end{document}